\documentclass[12pt]{article}
\usepackage[a4paper,left=2cm,right=1.2cm,top=3cm,bottom=4cm]
{geometry}
\usepackage{amsmath,amstext, amsthm, empheq, extpfeil,  
amsbsy,amssymb,marvosym,fancyhdr,graphicx,amscd,amsfonts,latexsym,delarray,stackrel,lineno,color,cite,appendix,xcolor,url,hyperref}
\usepackage{wasysym}
\usepackage{subfig}

\usepackage{srcltx}
\usepackage{mathrsfs}
\usepackage{float} 
\usepackage[all]{xy}
\usepackage{t1enc}
\usepackage{mathrsfs}
\usepackage{pifont}

\usepackage{mathpple}
\usepackage[T1]{fontenc}

\definecolor{Green}{rgb}{0,1,0}
\definecolor{Blue}{RGB}{0,0,191}
\definecolor{mathmodecolor}{RGB}{0,102,0}
\definecolor{keywordcolor}{RGB}{0,51,151}
\definecolor{sourcebackgroundcolor}{RGB}{255,247,223}
\definecolor{unixagred}{RGB}{255,0,0}
\definecolor{lightgray}{RGB}{191,191,191}
\definecolor{green}{RGB}{1,191,191}

\newcommand*\patchAmsMathEnvironmentForLineno[1]{%
  \expandafter\let\csname old#1\expandafter\endcsname\csname #1\endcsname
  \expandafter\let\csname oldend#1\expandafter\endcsname\csname end#1\endcsname
  \renewenvironment{#1}%
     {\linenomath\csname old#1\endcsname}%
     {\csname oldend#1\endcsname\endlinenomath}}%
\newcommand*\patchBothAmsMathEnvironmentsForLineno[1]{%
  \patchAmsMathEnvironmentForLineno{#1}%
  \patchAmsMathEnvironmentForLineno{#1*}}%
\AtBeginDocument{%
\patchBothAmsMathEnvironmentsForLineno{equation}%
\patchBothAmsMathEnvironmentsForLineno{align}%
\patchBothAmsMathEnvironmentsForLineno{flalign}%
\patchBothAmsMathEnvironmentsForLineno{alignat}%
\patchBothAmsMathEnvironmentsForLineno{gather}%
\patchBothAmsMathEnvironmentsForLineno{multline}%
}

\newtheorem{thm}{Theorem}[section]
\newtheorem{prop}[thm]{Proposition}
\newtheorem{cor}[thm]{Corollary}
\newtheorem{lem}[thm]{Lemma}

\newtheorem{fact}[thm]{Fact}

\def\Trace{{\rm Tr}}

\def\C{{\mathbb C}}

\def\F{{\mathbb F}}

\def\N{{\mathbb N}}

\def\Z{{\mathbb Z}}

\def\aarith{{\mathfrak A}}

\def\cO{{\mathcal O}}

\def\scal{{(\rnt,\cO)}}
\def\dar[#1]{\ar@<2pt>[#1]\ar@<-2pt>[#1]}
\def\qqq{\,,\,~\forall}

\newcommand{\ie}{{\it i.e.\/}\ }

\def\sin{{{\rm sin}}}
\def\cos{{{\rm cos}}}

\def\bm2{{\rm Bmod^2}}
\def\b2{{\rm Bmod^{\mathfrak s}}}

\newcommand{\nil}[1]{}

\def\rnt{{[0,\infty)\rtimes{\N^{\times}}}}

\def\aarith{{\mathscr A}}
\def\scal{{(\rnt,\cO)}}
\def\scal1{{\hat \aarith}}
\def\scal2{{\mathscr S}}

\title
{On an idea of Michael Atiyah}

\author{Alain Connes}

\begin{document}

\maketitle

\vspace{0.5cm}
\centerline{\em{In memoriam Michael Atiyah, }}

\centerline{\em{with admiration and gratitude}}
\vspace{0.3cm}

\begin{quote}
{ \em "In the broad light of day mathematicians check their equations and their proofs, leaving no stone unturned in their search for rigour. But, at night, under the full moon, they dream, they float among the stars and wonder at the miracle of the heavens. They are inspired. Without dreams there is no art, no mathematics, no life." }

(Michael Atiyah, Les Déchiffreurs 2008, Notices of the AMS, 2010).
\end{quote}

\section{Introduction}
The Feit--Thompson theorem on the solvability of finite groups of odd order was very much on Michael Atiyah's mind during his participation in the 2017 Shanghai conference on noncommutative geometry. Michael's lively presence there, and his inexhaustible enthusiasm for all mathematics --old, new
 and yet to be created-- were highlights of the meeting.
 
The idea of Michael's that we shall examine in this note was conceived by him during his flight home from Shanghai. It is a new strategy for FT, based on the iterative process sketched by him 
\begin{figure}[H]
\begin{center}
\includegraphics[scale=0.49]{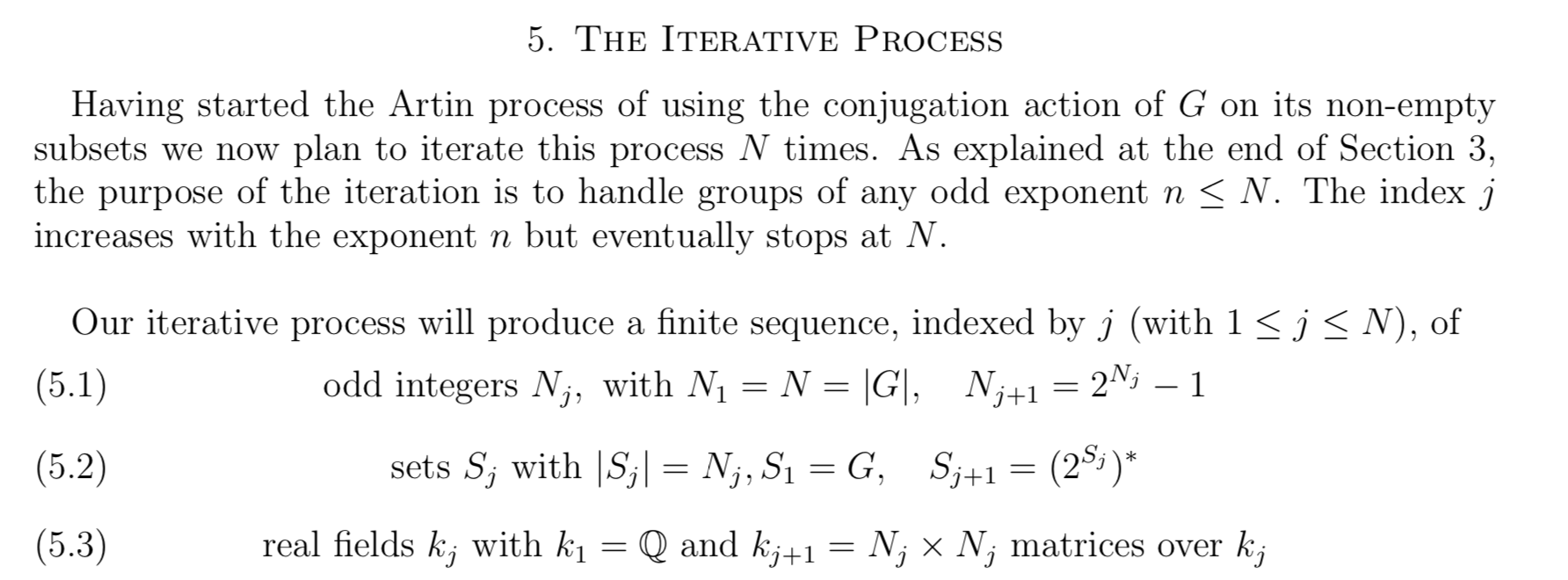}
\end{center}
\caption{ Extract from Michael Atiyah's notes \cite{Atiyah} \label{michaelnote} }
\end{figure} 
His idea is to use this process to construct a non-trivial character of the finite group $G$. Taken too literally this cannot work because it would apply to the 
group of permutations of $G$ which commute with the involution $g\mapsto g^{-1}$ and this group only has characters of even order.

The goal of the present paper, as a tribute to a luminous mathematical imagination that never dimmed, is to take seriously  his proposal and to show that, understanding it in a broader sense, one arrives at a very interesting idea. 

 My point of departure is to study more generally the iterations of the transformation which replaces a representation $\pi$ of a finite group $G$ on a finite dimensional complex vector space $E$ by the difference $\wedge \pi-1$, where $\wedge \pi$ is the action  of $G$ on the sum of exterior powers $\wedge^n E$ and $1$ is the trivial representation.

 We first show in \S \ref{oddremark} (Proposition \ref{wedgeG}) that the operation $\pi \mapsto \wedge \pi$ extends to virtual representations if and only if the group $G$ has odd order. We discuss briefly in \S \ref{euler} what goes wrong when the finite group has even order. In the odd case we show in Lemma \ref{adams} that  the wedge operation $\wedge: R(G)_\Z\to R(G)_\C$ is given by the formula:
 \begin{equation}\label{wedgeformulaintro}
 \wedge(x)=e^{\sum_1^k c_j\psi_j(x)}, \  \	c_j=(-1)^{j+1}\frac{H(j/k)}{k}
\end{equation}
where $k$ is the odd order of $G$, the $\psi_j$ are the Adams operations and the function $H(u)$ is the alternate sum  $\sum\frac{(-1)^n}{n+u}$. This allows one to extend  the wedge operation to complex valued class functions on the group $G$ and thus to make sense of the operation $\Psi(f):=\wedge f-1$ on such functions. We can then test Atiyah's idea of using iteration of $\Psi$.\newline
 Obviously the one-dimensional characters, \ie the morphisms $\chi:G\to \C^*$, are fixed points of $\Psi$ since $\wedge \chi=1+\chi$. We investigate the transformation $\Psi$ in the simplest non-trivial case, the finite odd non-abelian group of smallest order, which is of order $21$.  We first determine in \S \ref{smallest} the linear transformation $T$ on class functions for which the wedge is the composition of the exponential with $T$. We find that $T$ has a non-trivial kernel. In \S \ref{iterat}, we can finally investigate 
the iterations of the map $\Psi$ in our concrete example. As expected when one takes as a starting point the representation $\rho$ given by the action of $G$ on itself by conjugation, the iteration $\Psi^{\circ n}(\rho)$ blows up extremely fast and does not converge anywhere, thus making it hard to imagine that such  iteration can be used to prove the existence of a non-trivial one dimensional character. The nice surprise is that this divergent behavior is due to the bad choice of the initial data for the iteration. We show 
in \S \ref{iterat}, by a concrete computation, that for a whole class of initial data the iterations of the map $\Psi$ actually converge nicely to one of the two non-trivial characters of the group $G$. In more detail the group $G$ has 4 non-trivial conjugacy classes usually denoted as $7A$, $7B$ which are of order $7$ and size $3$ and $3A$, $3B$ which are of order $3$ and size $7$. We test first the iteration of $\Psi$ on central functions $f$ which are $1$ except on the classes $3A$, $3B$ on which they take conjugate values. Thus $f$ is complex valued and fulfills the condition $f(g^{-1})=\overline{f(g)}$. Such functions depend on a single complex number $z=f(3A)$ and the analysis of $\Psi$ becomes the iteration of a transformation of the complex plane which is given explicitly in Lemma \ref{explicitpsi}. 
\begin{figure}[H]
\begin{center}
\includegraphics[scale=0.8]{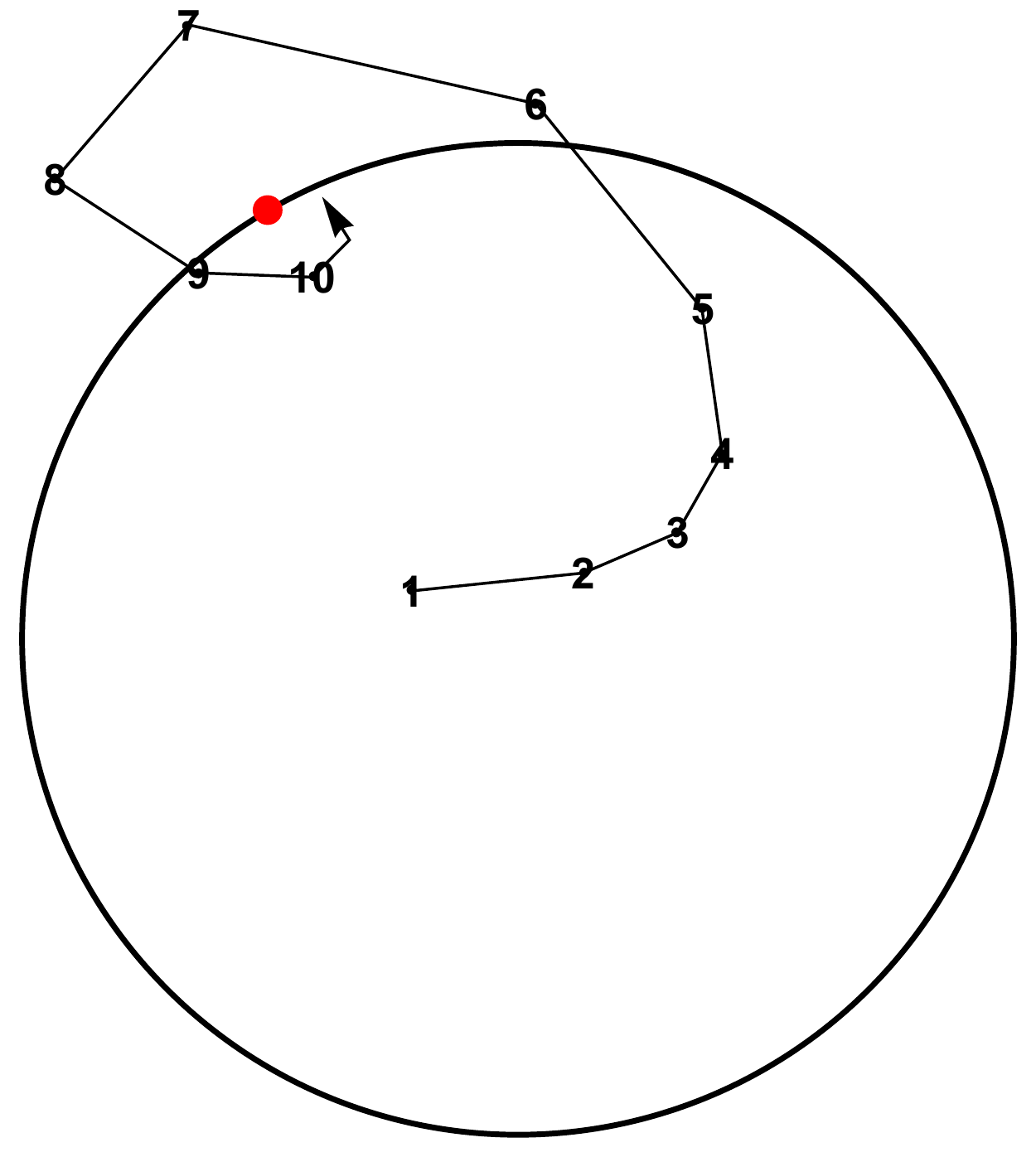}
\end{center}
\caption{ Iteration of the transformation $\Psi=$ wedge minus one. \label{matrixx} }
\end{figure} 
Concrete iteration of this transformation with initial data near the origin shows that starting from a point on the real axis one converges to the character of the trivial representation, but starting from a point in the upper half plane one converges to the non-trivial one-dimensional character of $G$ (and to its conjugate if one starts from a point in the lower half plane).
However, by investigating further the behavior of the transformation $\Psi$ on the non-trivial functions on the other conjugacy classes, we   find an attracting fixed point which does not correspond to a one dimensional character. We show in \S \ref{notall} that the fixed points for the other conjugacy class are solutions of a Lambert equation $we^w=u$ where the value $u=-\frac 67 2^{-\frac 57}\log 2\sim -0.362124$ happens to be in the interval $[-1/e,0)$ where two real solutions exist. While one of them corresponds to the character of the trivial representation, the other one does not come from a character. Moreover it is the fixed point that does not correspond to a character which is an attractor and this makes it quite difficult to conjecture a precise relation between fixed points and characters of one dimensional representations. 

Incidentally it happens interestingly that both the iteration of the exponentiation and the Lambert equation were considered by L. Euler \cite{euler1,euler2} and one may consider the determination of the fixed points of the operation $\Psi=$ wedge minus one, for finite groups of odd order as a variation on these early developments.

\tableofcontents 

\section{Oddness and the wedge}\label{oddremark}

We characterize odd finite groups \ie groups whose order is an odd integer, in terms of the $\lambda$-ring structure of the representation ring $R(G)_\Z$. The operation $\pi\mapsto \wedge \pi$ makes sense for finite dimensional representations and the issue is wether this map extends to virtual representation respecting the law $\wedge (\pi_1\oplus \pi_2)=\wedge(\pi_1)\otimes \wedge(\pi_2)$.
\begin{prop} \label{wedgeG}
	Let $G$ be a finite group, the following conditions are equivalent
 \begin{enumerate}
 \item The order of $G$ is odd.
 \item The map $\wedge$ extends to virtual representations as a morphism from the additive group $R(G)_\Z$ to the multiplicative monoid $R(G)_\C$. \end{enumerate}
\end{prop}
\proof Assume first that the map $\wedge$ extends to a morphism. One has $\wedge(0)=1$ and thus 
$$
\wedge(\pi)\wedge(-\pi)=1\in R(G)_\C
$$
Passing to the ring $C(G)$ of central functions by the Brauer map one gets that the character of 
$\wedge(\pi)$ is invertible in $C(G)$ and hence does not vanish. Thus by Lemma \ref{nochar} (below) it follows that $G$ has odd order. Assume now that $G$ has odd order. 
By construction the additive group $R(G)_\Z$ is obtained by symmetrization of the additive monoid of finite dimensional representations of $G$. Thus the existence of the extension is automatic provided the map $\wedge$ lands in invertible elements of the multiplicative monoid $R(G)_\C$. This monoid is the same as the monoid of complex valued class functions $C(G)$ under the pointwise product. An element in $C(G)$ is invertible iff it does not vanish and thus one gets the existence of the extension using Lemma \ref{nochar}.\endproof 

 \begin{lem}\label{nochar}
 Let $G$ be a finite group, the following conditions are equivalent
 \begin{enumerate}
 \item The order of $G$ is odd.
 \item For any finite dimensional representation $\pi$ of $G$ the character of $\wedge \pi$ does not vanish:
$$
\Trace(\wedge \pi(g))\neq 0\qqq g\in G
$$
 \end{enumerate}
 \end{lem}
\proof If $G$ has even order it contains an element $g$ of order $2$ and $\lambda(g)$ has the eigenvalue $-1$ in the regular representation, so that 
$$
\Trace(\wedge \pi(g))=\prod (1+\alpha)=0
$$ 
If $G$ has odd order any element $g\in G$ has odd order and hence $-1$ cannot be an eigenvalue of $\pi(g)$ since $(-1)^k\neq 1$ for $k$ odd. Thus
$$
\Trace(\wedge \pi(g))=\prod (1+\alpha)\neq 0
$$
Thus we get the desired conclusion.
\endproof

\section{Euler and divergent series}\label{euler}
We explain briefly what goes wrong in the even case. It is interesting to note that extending the map $\wedge$ hinges on the summation of divergent series. To be more precise one has by construction
\begin{equation}\label{wedgedefn}
\wedge= \sum_0^\infty \lambda_j	
\end{equation}
where the $\lambda_j$ are the lambda operations which are part of the structure of the representation ring. In the simplest case \ie when the group $G$ is reduced to one element, the operations $\lambda_j$ are given by the binomial  coefficients 
\begin{equation}\label{wedgedefn1}
 \lambda_j(x):=\frac{1}{j!}\prod_{0}^{j-1} (x-k)	
\end{equation}
Thus when we apply \eqref{wedgedefn} to negative numbers $x\in \Z$ we get the infinite sums 
$$
\wedge(-1)=1+(-1)+1+(-1)+1+\ldots 
$$
$$
\wedge(-2)=1-2 +3-4+5+\ldots 
$$
for which Leibniz and Euler give the values $\frac 12$ and $\frac 14$ respectively (see \cite{euler0,barbeau}).
\begin{figure}[H]
\begin{center}
\includegraphics[scale=0.6]{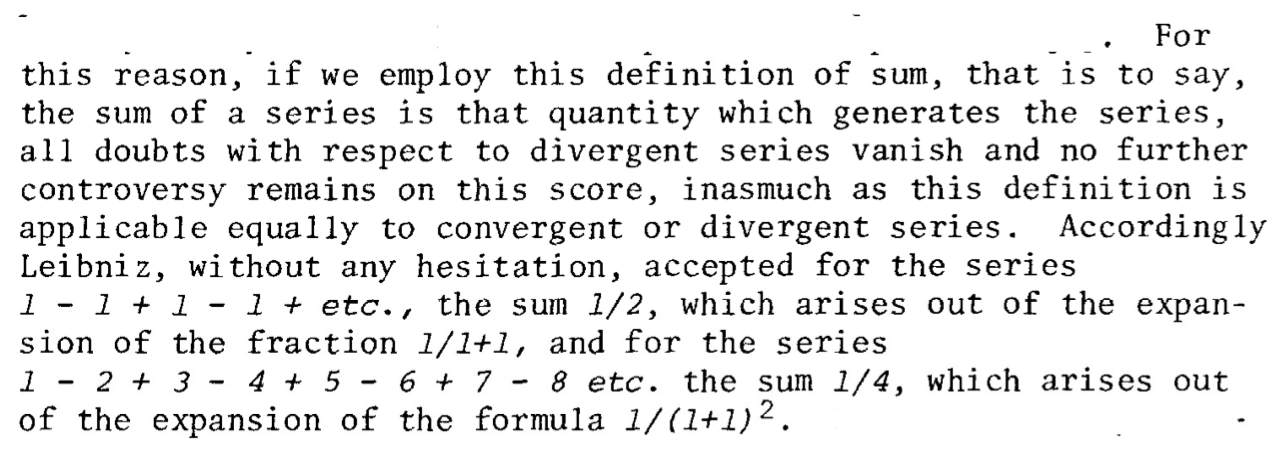}
\end{center}
\caption{Euler and divergent series. \label{graphofh} }
\end{figure}
These values give the correct answer $\wedge(-j)=2^{-j}$ thanks to the alternate signs. What Proposition \ref{wedgeG} shows is that one cannot define these sums of infinite series when all signs are the same  since the extension of the map $\wedge$  is not possible unless the group is odd. To see concretely what happens in the even case consider the cyclic group $C$ of order $2$. The elements of $R(G)_\Z$ are of the form $n+m\chi$ where $\chi^2=1$ and $n,m\in \Z$. The structure of $\lambda$-ring is such that 
\begin{equation}\label{wedgedefn2}
 \lambda_j(m \chi):=\frac{1}{j!}\prod_{0}^{j-1} (m-k)\chi^j	
\end{equation}
while the Brauer map sends in particular $\chi\mapsto -1$ by evaluation on the non-trivial conjugacy class. Thus on this class one gets the infinite sums
$$
\wedge(-\chi)\mapsto 1+1+1+\ldots 
$$
$$
\wedge(-2 \chi)=1+2 +3+4+5+\ldots 
$$
and Euler states clearly that such sums are infinite, unlike their alternate analogues. 
\section{Wedge in terms of Adams operations}\label{adams}

When the order of $G$ is odd it is possible to write a formula for the operation $\wedge$ as an exponential of a finite linear combination of the Adams operations in the representation ring $R(G)_\Z$. 
We use, for $u>0$ the following notation for the Hurwitz function
$$
H(u)=\sum _{m=0}^{\infty } \frac{(-1)^m}{m+u}
$$
where the sum is convergent and can be computed for rational values of $u$, such as 
$$
H(2/3)=\frac{\pi }{\sqrt{3}}-\log (2)
$$
In general one has with $k$
 an odd integer and $j\in \{1,\ldots ,k\}$,
 $$
   H(j/k)=(-1)^{j+1}\sum_{\vert v\vert<k/2} \left(\cos\left(\frac{2\pi jv}{k}\right) \log \left(2\cos(\frac{\pi v}{k}\right)+\frac{\pi v}{k}\sin\left(\frac{2\pi jv}{k}\right)\right)
   $$
 \begin{thm}\label{adams}
 Let $G$ be a finite group of odd order $k$. Let $\psi_j$ be the Adams endomorphisms of the representation ring $R(G)_\Z$. Then the wedge operation $\wedge: R(G)_\Z\to R(G)_\C$ is given by the formula:
 \begin{equation}\label{wedgeformula}
 \wedge(x)=e^{\sum_1^k c_j\psi_j(x)}, \  \	c_j=(-1)^{j+1}\frac{H(j/k)}{k}
\end{equation}
   \end{thm}
   \proof Both sides define morphisms from the additive group $R(G)_\Z$ to the multiplicative monoid $R(G)_\C$. Thus it is enough to check the equality \eqref{wedgeformula} when $x=\pi$ is a true representation of $G$. We evaluate the class function associated to both sides on a group element $g\in G$. By expressing the trace in terms of the eigenvalues $\alpha_j$ of $\pi(g)$, 
   the left hand side gives $\prod (1+\alpha_j)$. For the right hand side the formula is the exponential computed pointwise in the class ring $C(G)$  of 
   the expression obtained by replacing $\psi_j(x)$ by $\sum \alpha_i^j$. Thus since the exponential transforms addition into product, it is enough to show that for any $\alpha\in \C$ with $\alpha^k=1$ one has 
   $$
   1+\alpha=e^{\sum_1^k c_j\alpha^j}
   $$
   To see this we consider for $0<t<1$ the formula
   $$
  \log(1+\alpha t)= \sum_1^\infty (-1)^{j+1}\frac{\alpha^j t^j}{j}
   $$
   Using $\alpha^k=1$ we rewrite the right hand side as 
   $$
   \sum_1^k (-1)^{j+1}\alpha^j t^j \left(\sum_0^\infty (-1)^{u}\frac{ t^{ku}}{j+ku}\right)
   $$
   Applying the exponential to both sides (still for $t<1$) one has 
   $$
   1+t\alpha =Exp\left(\sum_1^k (-1)^{j+1}\alpha^j t^j \left(\sum_0^\infty (-1)^{u}\frac{ t^{ku}}{j+ku}\right)\right)
   $$
   Now when $t\to 1$ one has the convergence
   $$
   \lim_{t\to 1} \left(\sum_0^\infty (-1)^{u}\frac{ t^{ku}}{j+ku}\right)=\frac 1k H(j/k).
   $$
   Indeed one has $j\geq 1$ and 
   $$
   \frac{ t^{2kv}}{j+2kv}-\frac{ t^{2kv+k}}{j+2kv+k}=\frac{t^{2 k v} k}{(j+2kv)(j+2kv+k)}+\frac{t^{2 k v} \left(1-t^k\right)}{(j+2kv+k)}
   $$
   where one controls the behavior of the series involving the last term when $t\to 1$ because one has  $(1-t^k)\vert \log(1-t^{2k})\vert\to 0$ when $t\to 1$.
   Thus the continuity of the exponential allows one to pass to the limit when $t\to 1$ and gives the required equality. \endproof 
   \begin{cor} Let $G$ be a finite group of odd order $k$. The canonical extension of the wedge operation $\wedge: R(G)_\C\to R(G)_\C$ (given by the formula \eqref{wedgeformula}) 
 has a non-trivial kernel given by the central functions $f(g)$ such that
 $$
 \sum_1^k (-1)^{j}H(j/k)f(g^j)\in 2\pi i k\Z\qqq g\in G.
 $$  \end{cor}
 \proof The condition follows from \eqref{wedgeformula}. The transformation 
 \begin{equation}\label{linearma}
 T(f)(g):= \sum_1^k (-1)^{j}H(j/k)f(g^j)	
 \end{equation}
 is a linear map $T:R(G)_\C\to R(G)_\C$ on central functions. If the kernel of $T$ is trivial then $T$ is surjective and the inverse image of the lattice $(2\pi i k\Z)^c$, where $c$ is the number of classes, is non-trivial. Thus in all cases the kernel of $\wedge: R(G)_\C\to R(G)_\C$ is non-trivial.\endproof 
 \section{The first non-abelian odd group}\label{smallest}
One may wonder if the linear map \eqref{linearma} is always injective. This is not the case in general. We study the simplest non-abelian example. The first non-abelian group of odd order is of order $21$ and is the semi-direct product of the cyclic group $C(7)$ by the cyclic group $C(3)$ acting non-trivially. This is a subgroup of index $2$ in the affine group of $\F_7$. It is a transitive subgroup of the permutation group on $7$ letters and is denoted TransitiveGroup(7,3). 
This group has 5 conjugacy classes, the 4 non-trivial ones, labelled $\gamma_2,\gamma_3,\gamma_4,\gamma_5$ (or simply $2,3,4,5$ when no confusion can arise) are of the form:
$$
\gamma_2=7A, \gamma_3=7B,\ \gamma_4=3A,\ \gamma_5=3B,
$$
where the classes $7A$ and $7B$ are of order $7$ and size $3$, and the classes  $3A$ and $3B$ are of order $3$ and size $7$. The action of the inversion $j:g\mapsto g^{-1}$ simply interchanges $3A$ with $3B$ and $7A$ with $7B$. The action of $G$ on itself by conjugacy commutes with inversion and reduces (besides the trivial representation) to the sum of two copies of the representations of $G$  by permutations of sets with $3$ and with $7$ elements. The representation by permutations on the set with three elements $7A$ is by cyclic permutations using the morphism
$$
G=C(7)\rtimes C(3)\to C(3)
$$
The representation by permutations on the set with seven elements $3A$ is faithful and gives
$$
G\simeq {\rm TransitiveGroup}(7,3)
$$
(in the standard group theory notation).

The action of the $21$ power maps $g\mapsto g^n$, $n=1,\ldots ,21$, on the conjugacy class $2$ is given by
$$
\{2,2,3,2,3,3,1,2,2,3,2,3,3,1,2,2,3,2,3,3,1\}
$$
On the conjugacy class $3$ it is given by 
$$
\{3, 3, 2, 3, 2, 2, 1, 3, 3, 2, 3, 2, 2, 1, 3, 3, 2, 3, 2, 2, 1\}
$$
These classes are of order $7$ and thus the relevant exponent in \eqref{wedgeformula} reduces on the class $2$ to 
$$
\frac 17u(1)H(1)+\frac 17 u(2)\left(H(\frac 17)-H(\frac 27)-H(\frac 47)\right)
+\frac 17 u(3)\left(H(\frac 37)+H(\frac 57)-H(\frac 67)\right)
$$
and on the class $3$ to the same expression with $u(2)$ and $u(3)$ interchanged. But one has
$$
-H(\frac 17)+H(\frac 27)+H(\frac 47)=-H(\frac 37)-H(\frac 57)+H(\frac 67)
$$
In fact more precisely one has 
$$
H(\frac 17)+H(\frac 67)=\frac{\pi}{\sin(\frac \pi 7)},
H(\frac 27)+H(\frac 57)=\frac{\pi}{\cos(\frac{3 \pi }{14})},
H(\frac 37)+H(\frac 47)= \frac{ \pi}{ \cos \left(\frac{\pi }{14}\right)}
$$
and the above  equality follows from 
$$
\frac{1}{\sin(\frac \pi 7)}=\frac{1}{\cos(\frac{3 \pi }{14})}+\frac{ 1}{ \cos \left(\frac{\pi }{14}\right)}
$$
In fact one finds that in both cases the sum reduces to 
$$
\frac {\log 2}{7} \,  u(1)+\frac 37\log 2\left(u(2)+u(3)\right)
$$
Thus this means that the  $\wedge(u)$ only depends upon the sum $u(2)+u(3)$ of the central function $u$, \ie the sum of their values on the classes $7A$ and $7B$. The character table for $G$ is given by 

\vspace{0.1cm}

\begin{center}
\begin{tabular}{|l|c c c c r|}
\hline
  class & 1& 7A& 7B &3A & 3B \\
  \hline
  size & 1&3 & 3&7&7 \\
  \hline
  $\rho_1$ &  1&1&1 & 1& 1 \\
  $\rho_2$ &  1&1&1 & $j^2$ & $j$ \\
  $\rho_3$ &  1&1& 1 & $j$ & $j^2$  \\
  $\rho_4$ & 3&$\frac{-1-\sqrt{-7}}{2}$ &$\frac{-1+\sqrt{-7}}{2}$ & 0&0 \\
  $\rho_5$ & 3&$\frac{-1+\sqrt{-7}}{2}$&$\frac{-1-\sqrt{-7}}{2}$ &  0&0 \\
    \hline
\end{tabular}
\end{center}

\vspace{0.1cm}

where $j:=e^{\frac{2\pi i}{3}}$. One sees on this table that the characters of the representations $\rho_4$ and $\rho_5$ have the same sum on the classes $7A$ and $7B$. Thus we get:
\begin{lem}
The irreducible representations $\rho_4$ and $\rho_5$	of the group $G=C(7)\rtimes C(3)$ become isomorphic after applying the wedge, $\wedge \rho_4=\wedge \rho_5$.
\end{lem}
\proof Indeed the map $T$ of \eqref{linearma} takes the same value on $\rho_4$ and $\rho_5$ and the character of $\wedge \rho_j$ is a function of $T \rho_j$ so that the $\wedge \rho_j$ have the same character and hence are equal. \endproof 
In fact let us compute the characters of $\wedge \rho_j$. On the classes $7A$ and $7B$ they verify $u(2)+u(3)=-1$ thus one gets, since $u(1)=3$,
$$
\frac {\log 2}{7} \,  u(1)+\frac 37\log 2\left(u(2)+u(3)\right)=0
$$
After taking the exponential this gives the same contribution as $2\rho_1+\rho_4+\rho_5$ since the latter evaluates to $2-1=1$ on the classes $7A$ and $7B$. The answers agree for the trivial conjugacy class where they both give $8$. It remains to see what they give on the classes  $3A$ and $3B$. These classes are of order $3$ and the action of the power map on the class $4$ is of the form
$$
\{4, 5, 1, 4, 5, 1, 4, 5, 1, 4, 5, 1, 4, 5, 1, 4, 5, 1, 4, 5, 1\}
$$
 thus the relevant sum reduces on the class $4$ to 
$$
\frac 13u(1)H(1)+\frac 13 u(4)H(\frac 13)-\frac 13 u(5)H(\frac 23)
$$
and on the class $5$ to the same expression with $u(4)$ and $u(5)$ interchanged. But one has
$$
H(1/3)=\frac{\pi }{\sqrt{3}}+\log (2), \ H(2/3)=\frac{\pi }{\sqrt{3}}-\log (2)
$$
For both $\rho_4$ and $\rho_5$ one has $u(1)=3$, $u(4)=0$ and $u(5)=0$. Thus in both cases one gets 
$$
\frac 13u(1)H(1)+\frac 13 u(4)H(\frac 13)-\frac 13 u(5)H(\frac 23)=\frac 13u(1)H(1)=\log 2
$$
and after taking the exponential this gives $2$. Thus we have shown 
 \begin{equation}\label{wedgecompute}
 \wedge \rho_j=2\rho_1 +\rho_4+\rho_5\qqq j=4,5
 \end{equation}
and this corresponds in each case to the decomposition of $\wedge$ as $\sum \wedge^j$.
We have thus computed the matrix $T$ acting on column vectors as 
$$
T=\left(
\begin{array}{ccccc}
 \log (2) & 0 & 0 & 0 & 0 \\
 \frac{\log (2)}{7} & \frac{3 \log (2)}{7} & \frac{3 \log (2)}{7} & 0 & 0 \\
 \frac{\log (2)}{7} & \frac{3 \log (2)}{7} & \frac{3 \log (2)}{7} & 0 & 0 \\
 \frac{\log (2)}{3} & 0 & 0 & \frac{1}{3} \left(\log (2)+\frac{\pi }{\sqrt{3}}\right) & \frac{1}{3} \left(\log (2)-\frac{\pi }{\sqrt{3}}\right) \\
 \frac{\log (2)}{3} & 0 & 0 & \frac{1}{3} \left(\log (2)-\frac{\pi }{\sqrt{3}}\right) & \frac{1}{3} \left(\log (2)+\frac{\pi }{\sqrt{3}}\right) \\
\end{array}
\right)
$$

Its eigenvalues are 
$$
\left\{\frac{2 \pi }{3 \sqrt{3}},\log (2),\frac{6 \log (2)}{7},\frac{2 \log (2)}{3},0\right\}
$$
Its eigenvectors are the lines of the matrix:
$$
\left(
\begin{array}{ccccc}
 0 & 0 & 0 & -1 & 1 \\
 1 & 1 & 1 & 1 & 1 \\
 0 & 1 & 1 & 0 & 0 \\
 0 & 0 & 0 & 1 & 1 \\
 0 & -1 & 1 & 0 & 0 \\
\end{array}
\right)
$$

   \section{Iteration of $\Psi=$ wedge $-1$}\label{iterat}
 
The above computation restricted to class functions which take the same values on the classes $2=7A$ and $3=7B$ and the same values on $4=3A$ and $5=3B$ gives the following.

   \begin{lem}\label{group21}
 Let $G=C(7)\rtimes C(3)$. Let $f\in R(G)$ be a class function which is real and invariant under $g\mapsto g^{-1}$. Then for the three values of $f$ on the trivial class $1$, the class $2=7A$ and the class $4=3A$  one has
 \begin{equation}\label{wedgeformula1}
 \wedge(f)(1)=2^{f(1)}
 ,\ \wedge(f)(2)=2^{\frac 17 f(1)+\frac 67 f(2)},\ \wedge(f)(4)=2^{\frac 13 f(1)+\frac 23 f(4)}
\end{equation}
   \end{lem}
   \proof The formula for $\wedge f(2)$ comes from the above expression
   $$
\frac {\log 2}{7} \,  u(1)+\frac 37\log 2\left(u(2)+u(3)\right)
$$
with $u(2)=u(3)=f(2)$. The formula for $\wedge f(4)$ comes from 
$$
\frac 13u(1)H(1)+\frac 13 u(4)H(\frac 13)-\frac 13 u(5)H(\frac 23)
$$
with $u(4)=u(5)=f(4)$ while 
$
H(\frac 13)-H(\frac 23)=2\log 2
$.\endproof 
For the action of the group on itself by conjugation, the character associates to $g\in G$ the order of its commutant. Here this gives $21$ for the class $1$, $7$ for the classes $2=7A$, $3=7B$, and $3$ for $4=3A$ and $5=3B$. When one begins to iterate the wedge $-1$ on this starting point it gets immediately out of hands since after two steps it looks like, for the pairs $(f(1), f(4))$:
$$
\{21,3\},\ \{2097151,511\},\ \{2.272148509580683\times 10^{631305},4.674093841761078\times 10^{210537}\}
$$
However we now explain that this happens because of the wrong choice of starting point. We drop the condition $u(4)=u(5)$ and impose the conditions $u(1)=u(2)=u(3)=1$ which are by construction stable under the transformation $h\mapsto \wedge h-1$. We let $u(4)=x+iy$ be a complex number and replace the condition $u(4)=u(5)$ by  $\overline{u(4)}=u(5)$. Then one has
\begin{lem}\label{explicitpsi}
	The condition $\overline{u(4)}=u(5)$ is invariant under $h\mapsto \wedge h-1$ and  in terms of the real variables $(x,y)$, $u(4)=x+iy$, the transformation $h\mapsto \wedge h-1$ is given by
$$
\Phi(x,y):=\left(2^{\frac{1}{3} (2 x+1)} \cos \left(\frac{2 \pi  y}{3 \sqrt{3}}\right)-1,2^{\frac{1}{3} (2 x+1)} \sin \left(\frac{2 \pi  y}{3 \sqrt{3}}\right)\right)
$$
\end{lem}
\proof One has by definition, using $u(1)=1$, that the components   for $v=\wedge u-1$ are
$$
v(4)=Exp\left(\frac 13 H(1)+\frac 13 u(4)H(\frac 13)-\frac 13 u(5)H(\frac 23)\right)-1
$$
$$
v(5)=Exp\left(\frac 13 H(1)+\frac 13 u(5)H(\frac 13)-\frac 13 u(4)H(\frac 23)\right)-1
$$
One has 
$$
H(1)=\log (2), \ H(1/3)=\frac{\pi }{\sqrt{3}}+\log (2), \ H(2/3)=\frac{\pi }{\sqrt{3}}-\log (2)
$$
thus the exponents give respectively $\frac{\log 2}{3} (2 x+1)+\frac{2 i\pi  y}{3 \sqrt{3}}$ and its complex conjugate. Thus one gets the required formula. \endproof 
The Jacobian of $\Phi$ is given by 
$$
J(x,y)=\frac{1}{3} 2^{\frac{2 (x+2)}{3}}\left(
\begin{array}{cc}
\log (2) \cos \left(\frac{2 \pi  y}{3 \sqrt{3}}\right)  & -\frac{\pi  }{\sqrt{3}} \sin \left(\frac{2 \pi  y}{3 \sqrt{3}}\right)\\
 \log (2) \sin \left(\frac{2 \pi  y}{3 \sqrt{3}}\right) & \frac{\pi  }{\sqrt{3}}\cos \left(\frac{2 \pi  y}{3 \sqrt{3}}\right) \\
\end{array}
\right)
$$
When evaluated at the fixed point $P=(-1/2,\sqrt 3/2)$ one gets 
$$
(J(P)^*J(P))^{\frac 12}=\frac{1}{3} \left(
\begin{array}{cc}
\log (2)   & 0\\
0 & \frac{\pi }{\sqrt{3}} \\
\end{array}
\right)
$$
and one has the polar decomposition
$$
J(P)=\left(
\begin{array}{cc}
 \frac{1}{2} & -\frac{\sqrt{3}}{2} \\
 \frac{\sqrt{3}}{2} & \frac{1}{2} \\
\end{array}
\right).(J(P)^*J(P))^{\frac 12}
$$
It follows since both eigenvalues of $(J(P)^*J(P))^{\frac 12}$ are $<1$ that $P$ is an attractive fixed point. 

\begin{figure}[H]
\begin{center}
\includegraphics[scale=0.7]{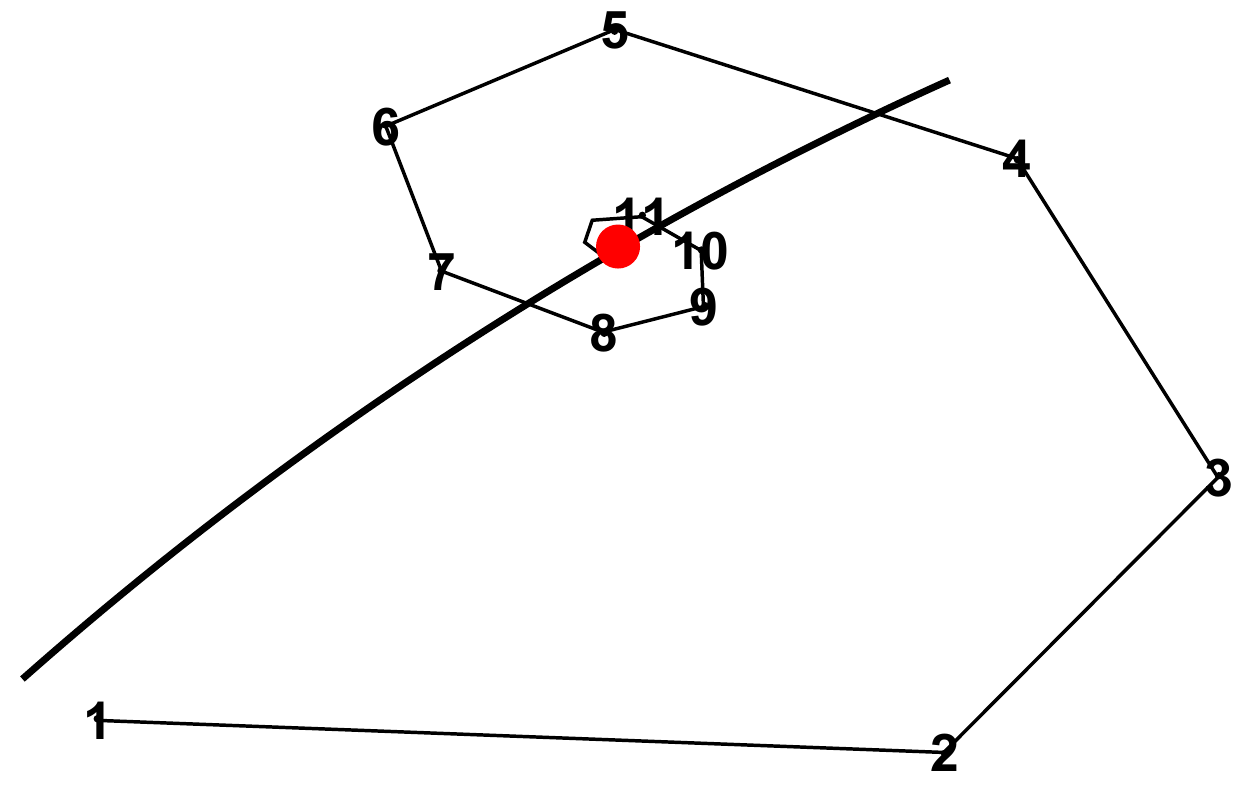}
\end{center}
\caption{ The transformation $\Psi$ spirals to the cubic root of $1$. \label{matrixx} }
\end{figure}

In fact the Jacobian of the transformation $\Psi$ at the fixed point $1\in \C$, is 
$$
\left(
\begin{array}{cc}
 \frac{4 \log (2)}{3} & 0 \\
 0 & \frac{4 \pi }{3 \sqrt{3}} \\
\end{array}
\right)
$$
one has $\frac{4 \log (2)}{3}<1$, $\frac{4 \pi }{3 \sqrt{3}}>1$ and this shows that as soon as one moves in the complex domain by a slight perturbation of the trivial representation one ends up in the domain of attraction of one of the two non-trivial characters.  Thus one has,
\begin{fact}
The iterates $\Psi^{\circ n}(1\pm i\epsilon)$ of a small complex  perturbation of the trivial representation on the classes 3A,3B, as above, converge, when $n\to \infty$ to a non-trivial character of $G$. 
\end{fact}

\section{Not all fixed points give characters}\label{notall}
\begin{figure}[H]
\begin{center}
\includegraphics[scale=0.7]{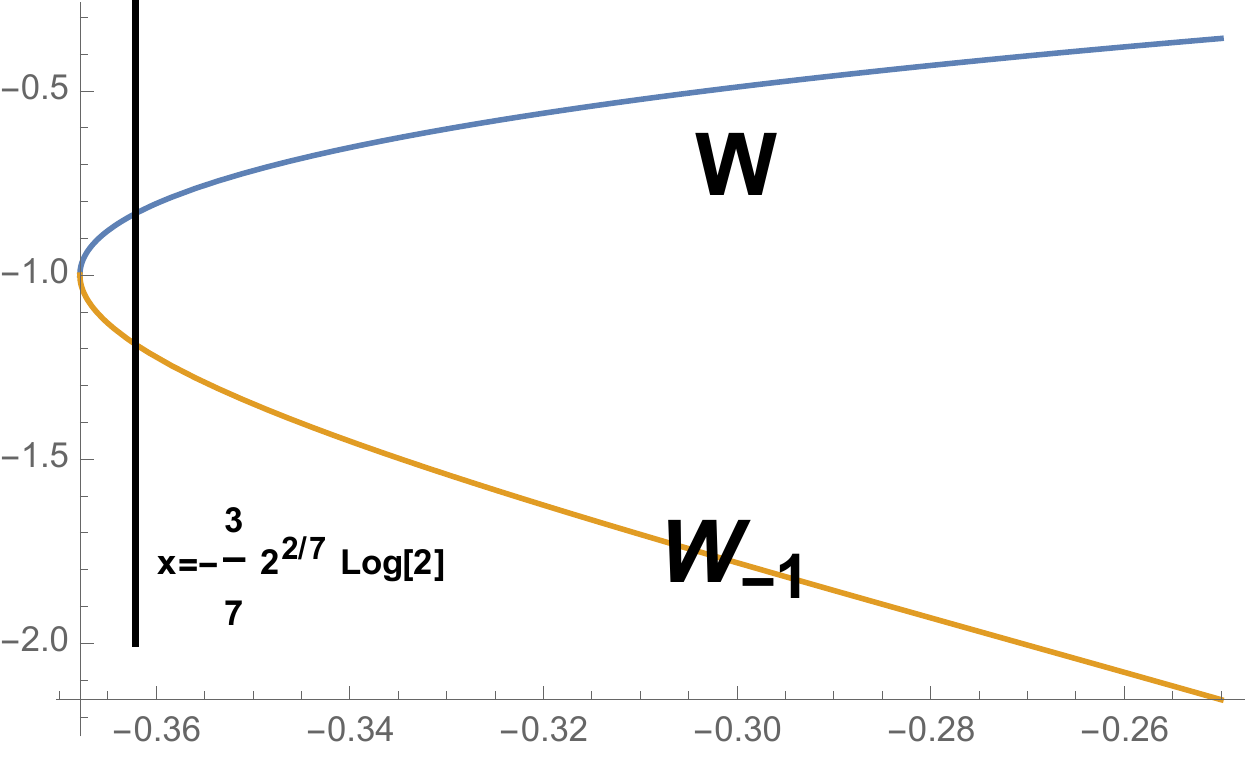}
\end{center}
\caption{ The two real branches of the Lambert function in the interval $[-1/e,0)$. \label{lambert} }
\end{figure}

This very nice behavior of the iteration of the map $\Psi$ is however not true in general and it will be enough to demonstrate that by looking at its behavior on the other two conjugacy classes $7A$ and $7B$. In that case as we saw above one gets after one iteration
the same value for $u(2)$ and $u(3)$ which is given, with $u(1)=1$, by 
$$
v(2)=v(3)=2^{\frac {1}{7} +\frac 37\left(u(2)+u(3)\right)}-1
$$ 
Thus in this case the transformation $\Psi=\wedge -1$ takes the form:
$$
\Psi(z)=2^{\frac {1}{7} +\frac 67 z}-1
$$
The fixed point equation $\Psi(z)=z$ is equivalent to $we^w=-\frac 67 2^{-\frac 57}\log 2$ where one lets $w=-\frac 67(1+z)\log 2$.  The numerical value $u=-\frac 67 2^{-\frac 57}\log 2\sim -0.362124$ is slightly larger than $-1/e\sim -0.367879$ and thus one falls in the range where two real solutions exist for the Lambert equation $we^w=u$ (see \cite{euler2,corless}). They are $W(u)$ and $W_{-1}(u)$ as shown in Figure \ref{lambert}.
 The obvious solution $z=1$ corresponds to $W_{-1}(u)$ but the main branch of the Lambert function $W(u)$ gives another non-trivial real solution $z\sim 0.401664$. This solution obviously does not correspond to a character and moreover it cannot be dismissed on the grounds of its nature as a fixed point of $\Psi$ because it turns out to be an attractor. Moreover the problem which arises is that the natural fixed point $1$ is not an attractor. In fact the derivative of $\Psi(z)$ at $z=1$ is $\frac{12 \log (2)}{7}$ which is $>1$. Thus $1$ is a repulsive fixed point while in fact the other one $z=-\frac{7 W\left(\frac{1}{7} (-3) 2^{2/7} \log (2)\right)}{6 \log (2)}-1$ which is $\sim 0.401664$ is an attractor. This fixed point gives a solution of the equation $\wedge f=1+f$ which even though it is real (\ie one has $f(g^{-1})=\overline{ f(g)}$ for all $g\in G$), and has dimension $1$ (\ie $f(1)=1$) does not correspond to a true representation of dimension $1$.

Note also that the Lambert equation $we^w=u$ admits an infinite number of complex solutions given by the various complex branches $W_k$, $k\notin \{-1,0\}$, of the Lambert function (see \cite{corless}). 
\section{Conclusion} 
From the above analysis we conclude that undoubtedly it is worth investigating in greater depth the fixed points of the map $\Psi:=\wedge -1$ on the complexified representation ring of finite groups of odd order. In particular the link with generalized forms of the Lambert equation should be explored.  But it remains totally unclear  whether this approach could lead to the construction of a non-trivial one dimensional character. One reason for skepticism is that only the power law on conjugacy classes is being used in deriving the formula for $\Psi$, while the compatibility with the power law does not suffice to select characters among unitary class functions as one sees in the simplest example of the odd abelian group $C(3)\times C(3)$.

\end{document}